\documentclass[12pt]{amsart}

\usepackage{amsfonts,amsmath,amssymb}

\usepackage{enumerate}

\numberwithin{equation}{section}

\newtheorem{theorem}{Theorem}[section]

\newtheorem{lemma}[theorem]{Lemma}

\newtheorem{proposition}[theorem]{Proposition}

\newtheorem{example}[theorem]{Example}

\newtheorem{definition}{Definition}[section]

\newtheorem{corollary}[theorem]{Corollary}

\newtheorem{remark}[theorem]{Remark}

\newcommand{\cl}[1]{\mathcal{#1}} 

\newcommand{\nor}[1]{\left\Vert #1\right\Vert}

\newcommand{\Map}[1]{\mathrm{Map}(#1)} 

\newcommand{\Lat}[1]{\mathrm{Lat}(#1)} 

\newcommand{\Alg}[1]{\mathrm{Alg}(#1)} 

\newcommand{\wsp}[1]{\overline{\mathrm{span}}^{\text{w*}}(#1)} 

\begin{document}

\title{Applications of operator space theory to nest algebra bimodules}

\author[G.~K.~Eleftherakis]{G. K. Eleftherakis}

\address{Department of Mathematics, University of Athens,
Panepistimioupolis 157 84, Athens, Greece.}

\email{gelefth@math.uoa.gr}



\date{}

\maketitle

\begin{abstract}Recently, Blecher and Kashyap have generalized the notion of $W^*$-modules over von Neumann algebras 
to the setting where the operator algebras are $\sigma -$weakly closed algebras of operators on a Hilbert space. They 
 call these  
modules {\em weak* rigged modules}. We characterize the weak* rigged modules over nest algebras. We prove that $Y$ 
is a right weak* rigged module over a nest algebra $\mathrm{Alg}(\cl M)$ if and only if there exists a completely 
isometric normal representation $\Phi $ of $Y$ and a nest algebra  $\mathrm{Alg}(\cl N)$ such that 
$\mathrm{Alg}(\cl N)\Phi (Y)\mathrm{Alg}(\cl M)\subset \Phi (Y)$ while $\Phi (Y)$ is implemented by a continuous nest 
homomorphism from $\cl M$ onto $\cl N.$  We describe some properties which are preserved by continuous CSL homomorphisms.  
\end{abstract}

\section{Introduction}

$C^*$-modules are well known as a generalization of the notion of Hilbert spaces introduced by Kaplansky \cite{kap}. They are 
very important tools for the study of $C^*$-algebras. The dual version is the notion of $W^*$-module. A $W^*$-module over a von 
Neumann algebra is a $C^*$-module which is a dual Banach space and defines a separately weak* continuous inner product \cite{bm}. 
 $W^*$-modules were originally introduced by Paschke \cite{pas}. They possess fruitful properties and 
appear in many versions such as weak* ternary rings of operators (TROs).

The development of operator space theory makes the generalization of $W^*$-modules possible even in the case of non-selfadjoint 
 dual operator algebras. In particular, Blecher and Kashyap provided a new characterization of $W^*$-modules 
which enabled the definition of an analogous notion for modules over non-selfadjoint 
 dual operator algebras \cite{bk}. They call these modules {\em weak* rigged modules}.

In the present paper we characterize the weak* rigged modules over nest algebras. We prove that a right dual operator 
module $Y$ over a nest algebra $B$ is a nest algebra bimodule, i.e., there exists a nest algebra $A$ such that 
$AYB\subset Y.$ Nest algebra bimodules are well constructed and described 
in \cite{ep}: If $ \mathrm{Alg}(\cl N),  \mathrm{Alg}(\cl M) $ are nest algebras and $Y$ is a weak* closed space satisfying 
 $\mathrm{Alg}(\cl N)Y \mathrm{Alg}(\cl M)\subset Y $, there exists an order preserving  left continuous map 
$\theta : \cl M\rightarrow \cl N$ such that $Y=Op(\theta ),$ where
\begin{equation}\label{op} Op(\theta )=\{y: \theta (M)^\bot yM=0\;\;\forall M\;\;\in \;\;\cl M\}.
\end{equation}
We prove that $Y$ is a right weak* rigged module over $\mathrm{Alg}(\cl M)$ if and only if there exists a normal completely isometric 
representation $\Phi $ and a nest algebra $\mathrm{Alg}(\cl N) $ such that $\Phi (Y)=Op(\theta )$ where 
$\theta : \cl M\rightarrow \cl N$ is a left-right continuous, therefore continuous, surjective nest homomorphism. 

We outline the results of this paper: In Section 2 we fix a continuous nest surjective homomorphism 
$\theta : \cl M\rightarrow \cl N.$ We prove that $Y=Op(\theta ),$ (see (\ref{op})) is a right weak* rigged module over 
 $\mathrm{Alg}(\cl M).$ There also exists an $\mathrm{Alg}(\cl N)-\mathrm{Alg}(\cl M) $ module $X$ such that $\mathrm{Alg}(\cl N)$ 
is isomorphic as an $\mathrm{Alg}(\cl N)$-module to the space $Y\otimes ^{\sigma h}_{\mathrm{Alg}(\cl M) }X$ which 
is the $\mathrm{Alg}(\cl M)$ balanced normal Haagerup tensor product of $Y$ and $X,$ \cite{elepaul}. Another result 
is that the algebra $\mathrm{Alg}(\cl N)$ is isomorphic to the algebra of the left multipliers over $Y, M_l(Y).$

In Section 3 we prove the converse: If $Y$ is an abstract right weak* rigged module over the nest algebra 
$ \mathrm{Alg}(\cl M) $, there exists a normal completely isometric representation $\Phi$ of $Y$, a nest algebra 
$\mathrm{Alg}(\cl N)$, and a continuous surjective nest homomorphism $\theta : \cl M\rightarrow \cl N$ 
with $\Phi (Y)=Op(\theta ). $

In Section 4 we present some examples and counterexamples of right weak* rigged modules over nest algebras. In the 
counterexamples we will see that the continuity of the nest is important.

In Section 5, inspired by the nest algebra case, we define the notion of a spatial embedding of a dual operator algebra 
$A$ in another dual operator algebra $B.$ In case $A$ and $B$ are CSL algebras which correspond to CSLs 
$\cl L_1, \cl L_2$, we prove that $A$ is spatially embedded in $B$ if and only if there exists a continuous CSL 
homomorphism from $\cl L_2$ onto $\cl L_1.$ Three natural consequences are the following:

(i) Let $\cl K(\mathrm{Alg}(\cl L_2) ),$ (resp. $\cl K(\mathrm{Alg}(\cl L_1)) $) be the subalgebra of 
$\mathrm{Alg}(\cl L_2)$ (resp. $\mathrm{Alg}(\cl L_1)$ )  which contains its compact operators. Then
$$ \mathrm{Alg}(\cl L_2)=\overline{\cl K(\mathrm{Alg}(\cl L_2))}^{w^*}\Rightarrow 
 \mathrm{Alg}(\cl L_1)=\overline{\cl K(\mathrm{Alg}(\cl L_1))}^{w^*}. $$

(ii) If $\mathrm{Alg}(\cl L_1)$ contains a non-zero compact operator (resp. finite rank operator), 
then $\mathrm{Alg}(\cl L_2)$ also contains  a non-zero compact operator (resp. finite rank operator).

(iii) If $\cl L_2$ is a synthetic lattice, then $\cl L_1$ is also a synthetic lattice.

\medskip
\medskip

In the following paragraphs we describe the notions we use in this paper; since we use extensively 
the basics of operator space theory, we refer the reader to the monographs \cite{bm}, 
\cite{er}, \cite{paul}, and \cite{pisier}  
for details. 

Let $H$ and $K$ be Hilbert spaces and $A\subset B(H)$ be an algebra. A subspace $X\subset B(K, H)$ 
is called a left module over $A$ if $AX\subset X.$ Similarly we can define 
right module over $A$. A left and right 
module over $A$ is called a bimodule over $A$. 
 An operator $Y$ is an abstract left (right) operator module over an abstract operator algebra $A$ if there exists  a 
 completely contractive bilinear map $A\times Y\rightarrow Y \;\;(Y\times A\rightarrow Y).$
 A left and right operator 
module over $A$ is called an operator bimodule over $A.$ 
 
  Let $A$ be a dual 
operator algebra and $Y$ be a dual operator space. We say that $Y$ is a left (right) 
dual operator module 
if the above completely contractive bilinear map is separately weak* continuous. A left and right dual operator 
module over $A$ is called a dual operator bimodule over $A.$ 

Two operator bimodules $Y$ and $Z$ over an operator algebra $A$ are 
called isomorphic as operator bimodules if there exists a completely 
isometric and surjective $A-$bimodule map $\pi : Y\rightarrow Z.$ This is then 
written $Y\cong Z$. If  $A$ is a dual operator algebra 
and $Y$ and $Z$ are dual operator $A-$bimodules, we write $Y\cong Z$ if also 
$\pi $ is weak* (bi)continuous.

If $Y$ is a dual right operator module over a dual operator algebra $B$ and $X$ is 
a left dual operator module over $B$, we denote by $Y\otimes ^{\sigma h}_B X$ the balanced 
normal Haagerup tensor product of $Y$ and $X$ which linearizes the separately 
weak* continuous completely bounded $B-$balanced 
bilinear maps \cite{elepaul}. If $Y$ (resp. $X$) is a left (resp. right) dual 
operator module 
over a dual operator algebra $A$, then $Y\otimes ^{\sigma h}_B X$ is also a left 
(resp. right) dual operator module over $A,$ \cite{elepaul}.

 The following is the definition of Morita equivalence used in this paper:

\begin{definition}\label{weak}\cite{bk} The dual operator algebras $B, A$ are called 
\textbf{weakly*  Morita equivalent} if there exist a $B-A$ dual operator module $X$ and 
an  $A-B$ dual operator module $Y$ such that $B\cong X\otimes ^{\sigma h}_AY$ and 
$A\cong Y\otimes ^{\sigma h}_BX$  as $B$ and $A$ dual operator bimodules respectively. 
\end{definition}

Fix $n, m\in \mathbb{N}.$ If $X$ is an operator space, $M_{m, n}(X)$ is the operator space of $m \times n$ 
matrices with entries in $X.$  Then $ M_n(X)=M_{n,n}(X), \;\;C_{n}(X)=
M_{n, 1}(X),\;\; R_{n}(X)=M_{1, n}(X).$ 

If $X$ is a subspace of $B(H, K)$, where $H$ and $K$ are Hilbert spaces, we denote by 
$R^{fin}_\infty (X)$ (resp. $C^{fin}_\infty (X)$) 
the space of operators $(x_1, x_2,...): H^\infty \rightarrow K$ (resp. $(x_1, x_2,...)^T: 
H \rightarrow K^\infty $) such that $x_i\in X$ for all $i$ and there exists $n_0\in 
\mathbb{N}$ for which $x_n=0$ for all $n\geq n_0.$

In this paper we shall use the notion of the multiplier algebra $M_l(Y)$ of a dual operator space $Y$ \cite[4.5.1]{bm}. We 
recall that $M_l(Y)$ is a unital dual operator algebra and $Y$ is a left dual operator module over $M_l(Y).$ 
Every $u\in M_l(Y)$ is a weak* continuous map from $Y$ into $Y.$ A bounded net $(u_t)_t\subset M_l(Y)$ 
converges in the weak* topology to $u\in M_l(Y)$ iff $u_t(y)\stackrel{w^*}{\rightarrow }u(y)$ for all $y\in Y.$ 

We present the definition of right weak* rigged modules over dual operator algebras. Throughout this paper we shall use 
this definition or other equivalent definitions from \cite{bk2}.

\begin{definition}\label{themelion} Suppose that $Y$ is a dual operator space and a right 
module over a dual operator algebra $B.$ Suppose that there exists a net of positive integers $(n(\alpha ))$ and 
weak* continuous completely contractive $B-$module maps $$\phi _\alpha : Y\longrightarrow C_{n(\alpha )}(B) ,\;\;\; 
\psi _{\alpha }: C_{n(\alpha )}(B) \rightarrow Y,$$ with $\psi_\alpha ( \phi_\alpha (y))\rightarrow y $ in the 
weak* topology on $Y,$ for all $y\in Y.$ Then we say that $Y$ is a \textbf{right weak* rigged module over $B$}.
\end{definition}

If $\cl L$ is a set of projections of a Hilbert space $H$ we denote by $\Alg{\cl L}$ the algebra 
 $$ \{x\in B(H): N^\bot xN=0\;\;\forall \;\;N\;\;\in \cl L \}. $$
A \textbf{nest} $\cl N$ is a totally ordered set of projections of a Hilbert space 
$H$ containing the zero and identity operators which is closed under arbitrary 
intersections and closed spans. The corresponding \textbf{nest algebra} is 
$ \mathrm{Alg}(\cl N) .$
If $N\in \cl N$, we denote by $N_-$ the projection onto the closed span of the 
union $\cup _{\stackrel{M<N}{M\in \cl N}}(M(H)).$ If $N_-<N$ we call the projection 
$N\ominus N_-$ an \textbf{atom}. 
If a nest has no atoms, it is called a \textbf{continuous nest}. If the atoms span 
the identity operator, the nest is called a \textbf{totally atomic nest}. An order preserving 
map  
between two nests is called a \textbf{nest homomorphism}.  If this map is injective and surjective, 
 it is called a \textbf{nest isomorphism}.

If $\cl N_1$ and $\cl N_2$ are nests acting on the Hilbert spaces $H_1, H_2$ 
respectively, and if $\theta : \cl N_1\rightarrow \cl N_2$ is a nest homomorphism, 
we denote by $Op(\theta )$ the space of operators $x\in B(H_1, H_2)$ 
satisfying $\theta (N)^\bot x N=0$ for all $N\in \cl N_1.$ Observe that 
$Op(\theta )$  is an $ \mathrm{Alg}(\cl N_2)-\mathrm{Alg}(\cl N_1) $ bimodule. 
If $Y$ is a subspace of $B(H_1, H_2),$ the map sending every 
projection $p$ of $H_1$ to the projection of $H_2$ generated by the vectors 
of the form $yp(\xi ): y\in Y, \xi \in H_1$ is denoted by $ \Map{Y} $.

Finally, if $X$ is a normed space, we denote by $Ball(X)$ the unit ball of $X.$

We recall the following results which will be used in this paper:

\begin{theorem}\label{A1} \cite{ep} Let $\cl N_1, \cl N_2$ be nests, $Y$ be a weak* closed 
 $ \mathrm{Alg}(\cl N_2) - \mathrm{Alg}(\cl N_1) $ bimodule, and $\theta$ be the restriction of $\Map{Y} $ 
to $\cl N_1.$ Then $\theta $ is left continuous, $\theta (\cl N_1)\subset \cl N_2$, and $Y=Op(\theta ).$ 
\end{theorem}

\begin{theorem}\label{A2} \cite{ele3} Let $\cl N_1$ and $ \cl N_2$ be nests, $B$ and $ A$ be the corresponding nest algebras, 
$\theta : \cl N_1\rightarrow \cl N_2$ be a nest isomorphism, and $Y=Op(\theta ), X=Op(\theta ^{-1})$. Then 

(i) $B=\wsp{XY}$ and $A=\wsp{YX}.$ 

(ii) $B$ and $A$ are weakly* Morita equivalent. In particular, $B\cong X\otimes ^{\sigma h}_AY$ as dual operator $B-$modules,
and $A\cong Y\otimes ^{\sigma h}_BX$ as dual operator $A-$modules.    
\end{theorem}

The last theorem implies, \cite[Theorem 3.3]{bk}, the following:

\begin{theorem}\label{A3} Let $A, B, X, Y$ be as in Theorem \ref{A2}. 
There exist nets $(y_t)_t\subset Ball(R_\infty ^{fin
}(Y)), $$(x_t)_t\subset Ball(C_\infty ^{fin
}(X)) $ such that $y_tx_t\stackrel{w^*}{\rightarrow }I_A$ where $I_A$ is the identity of $A,$ 
and nets $(u_i)_i\subset Ball(R_\infty ^{fin
}(X)), $$(w_i)_i\subset Ball(C_\infty ^{fin
}(Y)) $ such that $u_iw_i\stackrel{w^*}{\rightarrow }I_B$ where $I_B$ is the identity of $B.$ 

 \end{theorem}

\section{Continuous nest homomorphisms}

Let $\cl N, \cl M$ be nests acting on the Hilbert spaces $L, H$ respectively, $A$ and $B$ be the corresponding 
nest algebras, and $\phi : \cl M\rightarrow \cl N$ be a nest homomorphism. We assume that $\cl N$ is generated as a nest by $\phi (\cl M
)$. Observe that if $\phi $ is  continuous, then $\phi (\cl M)=\cl N.$ We define the following spaces:
$$Y=Op(\phi )=\{y\in B(H, L): \phi (M)^\bot yM=0, \;\;\forall \;\;M\in \cl M\}$$
$$X=\{x\in B(L, H): M^\bot x\phi (M)=0, \;\;\forall \;\;M\in \cl M\}. $$
 Observe that $Y$ is an $A-B$ bimodule and $X$ is a $B-A$ bimodule.

\begin{theorem}\label{B1} If $\phi $ is continuous, then

(i) $A=\wsp{YX}$, 

(ii) the space $Y$ is a right weak* rigged module over $B$,

(iii) the algebra $A$ is isomorphic as a dual $A-$operator module with $Y\otimes ^{\sigma h}_BX$,

(iv) the algebra $A$ is isomorphic as a dual operator algebra with the algebra $M_l(Y)$ of left 
multipliers over $Y.$ 
\end{theorem}
\textbf{Proof}

We define the nests
$$\cl L_1=\cl M^{(2)},\;\;\; \cl L_2=\{M\oplus \phi (M): M\in \cl M\}$$
 and the corresponding nest algebras

$$\Alg{\cl L_1}=\left(\begin{array}{clr} B &   B \\  B & B \end{array}\right)$$
$$\Alg{\cl L_2}=\left(\begin{array}{clr} B &   X \\  Y & A \end{array}\right). $$
The map $\psi : \cl L_1\rightarrow \cl L_2: M\oplus M\longrightarrow M\oplus \phi (M)$ is a nest 
isomorphism. We can easily check that 
$$\hat{Y}=Op(\psi )=\left(\begin{array}{clr} B &   B \\  Y & Y \end{array}\right)$$
 and 
$$\hat{X}=Op(\psi ^{-1})=\left(\begin{array}{clr} B &   X \\  B & X \end{array}\right). $$
Theorem \ref{A2} implies that $\Alg{\cl L_2}=\wsp{\hat{Y}\hat{X}}.$ Therefore $A=\wsp{YX}.$ 
Using Theorem \ref{A3} we find nets $$(\hat{y}_t)_t\subset Ball( R_\infty ^{fin
}(\hat{Y})), \;\;\;(\hat{x}_t)_t\subset Ball( C_\infty ^{fin
}(\hat{X})) $$ with the property $\hat{y}_t\hat{x}_t \stackrel{w^*}{\rightarrow }I_{\Alg{\cl L_2} }.$ 

We may take 
$$ \hat{y}_t =\left(\begin{array}{clr} b_t^1 & b_t^2 \\ y_t^1 & y_t^2 \end{array}\right) $$
where $ b_t^1,  b_t^2\in  R_\infty ^{fin}(B), \;\;\;y_t^1,  y_t^2 \in   R_\infty ^{fin}(Y) $ and  
$$\hat{x}_t=\left(\begin{array}{clr} c_t^1 & x_t^1 \\ c_t^2 & x_t^2 \end{array}\right) $$
where $ c_t^1,  c_t^2\in  C_\infty ^{fin}(B), \;\;\;x_t^1,  x_t^2 \in   C_\infty ^{fin}(X). $ 

Put $$ y_t=(y_t^1,  \;\;\;y_t^2)\in R_\infty ^{fin}(Y), \;\;\; x_t=(x_t^1, \;\;\;  x_t^2)
^T\in C_\infty ^{fin}(X). $$ 
Since $\hat{y}_t $ and $\hat{x}_t $ are contractions, $y_t$ and $x_t$ are also contractions. Clearly the net  
$$y_tx_t=y_t^1x_t^1+y_t^2x_t^2$$ converges to the identity of the algebra $A.$ We deduce that $Y$ is a right 
weak* rigged module over $B.$ (See the fourth description of the weak* rigged modules in \cite{bk2}).

Statements (iii) and (iv) follow from (ii) using the results of  \cite{bk2}. We include  a proof for completeness.

The map $$Y\times X\rightarrow A,\;\;\; (y,x)\rightarrow yx$$ is a complete contraction, weak* continuous, and  
an $A-$module map. 
So it induces a map $$\rho : Y\otimes ^{\sigma h}_BX \rightarrow A$$ which is also a complete contraction, 
weak* continuous, and an $A-$module map. We claim that $\rho $ is $1-1.$ Indeed, if $y\in Y$ and $ x\in X$, then 
$$(y\otimes _Bx)y_tx_t=y\otimes _B(xy_tx_t)=(yxy_t) \otimes_B x_t=\rho (y\otimes _Bx)( y_t\otimes_B x_t ). $$
 It follows that 
$$wy_tx_t=\rho (w)(y_t\otimes_B x_t )\;\;\;\forall w\in Y\otimes ^{\sigma h}_BX . $$
Since  $(y_tx_t)_t$ converges to the identity operator, $\rho (w)=0$ implies $w=0.$

If $v_1, v_2,...,v_m\in X, w_1,w_2,...,w_m\in Y$, we have 
\begin{align*}& \nor{\sum_{k=1}^mw_k\otimes _B(v_ky_tx_t)}=\nor{(\sum_{k=1}^mw_kv_ky_t)\otimes _Bx_t}\\ \leq &
\nor{\sum_{k=1}^mw_kv_ky_t}\|x_t\|\leq \nor{\sum_{k=1}^mw_kv_k}
\end{align*}
Letting $y_tx_t\rightarrow I_A$, we have 
$$\nor{\sum_{k=1}^mw_k\otimes _Bv_k}\leq \nor{\sum_{k=1}^mw_kv_k}. $$
We proved that the restriction of $\rho $ to $$ Y\otimes^h _BX =\overline{span}^{\|\cdot\|}\{
 y\otimes _Bx: y\in Y, x\in X \}$$ is an isometry. Choose $u\in  Y\otimes^{\sigma h }_BX $ and a net 
$(u_k)_k$ in $ Y\otimes^h _BX $ converging to  $u$ in  the weak* topology. For any finite rank operators 
$g, f\in Ball(A)$ we have 
$$\rho (fu_kg)=f\rho (u_k)g  \stackrel{\|\cdot\|}{\longrightarrow } f\rho (u)g=\rho (fug). $$
Since $\rho $ is $1-1$, we conclude that $fug\in Y\otimes^h _BX $, and so 
$$\|fug\|=\|\rho (fug)\|=\|f\rho (u)g\|\leq \|\rho (u)\|. $$ 
 The identity of $A$ is in the weak* closure of its finite rank contractions, \cite{dav}. Therefore 
$\|u\|\leq \|\rho (u)\|.$ 
We proved that $\rho $ is an isometry. Similarly we can prove that it is a complete isometry. 

It remains to prove statement (iv). Define $\sigma : A\rightarrow M_l(Y)$ by 
$$\sigma (a)(y)=ay\;\; \forall \;\;a\in A, y\in Y$$ which is clearly contraction. Let $u$ be in $M_l(Y).$ 
Since $u$ is a weak* continuous $B-$module homomorphism,we have
$$\lim_t\sigma (u(y_t)x_t)(y)=\lim_tu(y_t)x_ty=\lim_tu(y_tx_ty)=u(y),\;\;\forall \;\;y\in Y. $$
It follows that $(\sigma (u(y_t)x_t))_t$ converges in the weak* topology of $M_l(Y) $ to $u.$ So 
$u\in \overline{\sigma (A)}^{w^*}.$ If $a\in A$,
$$\|ay_tx_t\|=\| \sigma (a)(y_t)x_t\|\leq \|\sigma (a)\|. $$ The equality $\|a\|=\|\sigma (a)\|$ follows from the fact  
$ay_tx_t\stackrel{w^*}{\longrightarrow }a.$ So 
 $\sigma $ is an isometry. Similarly 
we can prove that it is a complete isometry. The Krein--Smulian Theorem, \cite[A.2.5]{bm}, implies 
that $$M_l(Y)=\overline{\sigma (A)}^{w^*}=\sigma (A). \qquad \Box$$

\begin{theorem} \label{B2} The following are equivalent:

(i) $\phi $ is continuous;

(ii) $A=\wsp{YX}$;

(iii) The identity of the Hilbert space $L$ belongs to $\wsp{YX}.$
\end{theorem} 
\textbf{Proof} Statement (i) implies (ii) by Theorem \ref{B1}. Since nest 
algebras are unital, (ii) implies (iii).  Suppose that $I_L\in \wsp{YX}.$ There 
exist nets $ (w_t)_t\subset R_\infty ^{fin}(Y), (u_t)_t\subset C_\infty ^{fin}(X) $ such that $w_tu_t\stackrel{w^*}{\longrightarrow 
}I_A.$ If $(M_i)_i$ is a decreasing net of projections in $\cl M$ converging in $M,$ we may assume that 
the net $(\phi (M_i))_i$ converges in $N\in \cl N.$ Clearly $\phi (M)\leq N.$ We shall prove that $N\leq \phi (M)$.

$$ N=\lim_tw_tu_tN=\lim_t(\lim_i w_tu_t\phi (M_i) )=\lim_t(\lim_iw_tM_i^\infty u_t\phi (M_i) ) . $$
Since $(M_i^\infty )_i$ (resp. $(\phi (M_i))_i)$ converges in the strong operator topology to 
$M^\infty $ (resp. to $N$),
\begin{align*} N=& \lim_tw_tM^\infty u_tN=\lim_t\phi (M) w_tM^\infty u_tN\\=&\phi (M)
\lim_t(w_tM^\infty u_tN )\Rightarrow N\leq \phi (M). 
\end{align*}   
We proved that $\phi(M)=N.$ Therefore $\phi $ is right continuous. The proof of the fact $\phi $ is left continuous is similar.$\qquad \Box$

\medskip

Let $ \stackrel{\sim }{X} $ be the space $\:_{w^*}CB(Y, B)_B.$ This  is the space  of weak* continuous completely 
bounded $B-$module maps from $Y$ to $B$ and it is, \cite{bk}, a dual $B-A$ module under the actions

$$B\times \stackrel{\sim }{X} \longrightarrow \stackrel{\sim }{X} : (b,u)\rightarrow b\cdot u,\;\;\;  b\cdot u(y)=bu(y)\;\;\forall 
y\in Y,  $$
$$ \stackrel{\sim }{X} \times A\longrightarrow \stackrel{\sim }{X} : (u,a)\rightarrow u\cdot a,\;\;\;  u\cdot a(y)=u(ay)\;\;\forall 
y\in Y.  $$

\begin{theorem} \label{B3} Suppose that $\phi $ is continuous. Then:

(i) The map $\iota : X\rightarrow \stackrel{\sim }{X} : \iota (x)(y)=xy$ is a weak* continuous completely isometric 
 surjective map.

(ii) $M^\bot \cdot \stackrel{\sim }{X} \cdot \phi (M)=0$ or equivalently $M^\bot u(\phi (M)y)=0$ for all 
$M\in \cl M,\;\; u\in \stackrel{\sim }{X},\;\; 
y\in Y.$ 
\end{theorem}
\textbf{Proof}
Clearly $\iota $ is a complete contraction. By Theorem \ref{B1} there exist nets $$ (y_t)_t\subset Ball(R_\infty ^{fin}(Y)),
(x_t)_t\subset Ball(C_\infty ^{fin}(X)) $$ such that $y_tx_t\stackrel{w^*}{\rightarrow }I_A.$ So 
for all $x\in X$ we have 
$$\|xy_tx_t\|=\|\iota (x)(y_t)x_t\|\leq \|\iota (x)\|. $$ Since $w^*-\lim_txy_tx_t=x$, we have 
$\|x\|\leq \|\iota(x)\|. $ Hence $\iota $ is an isometry. Similarly, we can prove that $\iota $ is a complete isometry. If $u\in 
 \stackrel{\sim }{X} $ for all $y\in Y$, 
\begin{align*} & u(y_tx_ty)=u(y_t)x_ty=u(y_t)\iota (x_t)(y) \\ \Rightarrow & u(y)=\lim_tu(y_t)\iota (x_t)(y)=\lim_t\iota (
u(y_t)x_t)(y).
\end{align*}
 Thus $$u=w^*-\lim_t\iota (u(y_tx_t))\in \overline{\iota (X)}^{w^*}$$ in the 
weak* topology of $CB(Y, B).$ By the Krein--Smulian Theorem 
$\iota (X)= \stackrel{\sim }{X},$ so statement (i) holds. Since $M^\bot X \phi (M)=0$ for all $M\in \cl M$, 
(i) implies (ii). $\qquad \Box$

\medskip

Theorem \ref{B6} describes how we can construct right weak* rigged modules over nest algebras. 
 We first need the following lemma:

\begin{lemma}\label{B7} Let $\cl M, \cl N$ be nests acting on Hilbert spaces $H, L$ respectively, and   
let $B=\Alg{\cl M}, A=\Alg{\cl N} $ be the corresponding nest algebras. Suppose that there exist spaces $X_0, Y_0$ such that:

(i) $X_0AY_0\subset B$;

(ii) The identity of the Hilbert space $L$ belongs to $\wsp{Y_0X_0}.$
 
Then $Y=\wsp{Y_0B}$ is a right weak* rigged module over $B.$ 
\end{lemma}
\textbf{Proof} Clearly $Y$ is  a left $B-$module. We shall prove that it is a right $A-$module:

\begin{equation}\label{mod} AY_0B\subset \wsp{Y_0X_0}AY_0B\subset \wsp{Y_0X_0AY_0B}.
\end{equation}

Since $X_0AY_0\subset B,$ (\ref{mod}) implies that 
$$AY_0B\subset Y\Rightarrow AY\subset Y. $$
For every $M\in \cl M$, we write $\phi (M)$ for the projection onto $\overline{Y_0M(H)}.$ Theorem \ref{A1} implies 
that $\phi (M)\in 
\cl N$ and $Y=Op(\phi ).$ 

Let $X$ be the space $$\{x\in B(L,H): M^\bot x\phi (M)=0\;\;\forall M\;\;\in\;\;\cl M \},$$ 
and $\hat{X}$ be the space $\wsp{BX_0A}.$ If $M\in \cl M$, 
$$ M^\bot BX_0A Y_0M\subset M^\bot BM=0\Rightarrow M^\bot BX_0A \phi (M)=0. $$
We proved that $\hat{X}\subset X.$ On the other hand,
$$\wsp{Y\hat{X}}=\wsp{YBX_0A}\supset \wsp{Y_0X_0}A\supset A. $$
It follows that $$\wsp{YX}\supset A\Rightarrow I_L\in \wsp{YX}. $$ Theorems \ref{B1} and \ref{B2} imply that 
$Y$ is a weak* rigged module over $B. \qquad \Box$ 
 
\begin{theorem}\label{B6}Let $\cl M$ be a nest acting on the Hilbert space $H$  and   
$B=\Alg{\cl M} $ be the corresponding nest algebra. Suppose that there exist spaces $X_0\subset B(L, H), Y_0\subset B(H, 
L)$ such that:

(i) $X_0Y_0\subset B$;

(ii) The identity of the Hilbert space $L$ belongs to $\wsp{Y_0X_0}.$
 
It follows that the space $Y=\wsp{Y_0B}$ is a right weak* rigged module over $B.$ 
 
\end{theorem}
\textbf{Proof} We define the map
$$\theta : \cl M\rightarrow pr(B(L)), \;\;\theta (M)=\overline{Y_0M(H)}$$
and the  algebra
$$\cl A=\Alg{\theta (\cl M)}=\{a:\;\;\; \theta(M)^\bot a \theta(M)=0, \;\;\forall M\;\;\in \;\;\cl M \}. $$
 Since  $\theta (\cl M)$ is a totally ordered set of projections, $A$ is a nest algebra. 
 By Lemma \ref{B7} it suffices to prove 
that $ X_0AY_0\subset B .$

For all $M\in \cl M, x\in X_0$ we have
$$x\theta (M)(H)\subset \overline{xY_0M(H)}\subset \overline{BM(H)}\subset M(H). $$
So $M^\bot X_0\theta (M)=0 \;\;\forall\;\; M\;\;\in\;\; \cl M.$
It follows that 
\begin{align*} & M^\bot X_0AY_0M = M^\bot X_0A\theta (M)Y_0M\\= & M^\bot X_0\theta (M)AY_0M =0 ,\;\;
\forall  \;\;M\;\;\in \;\;\cl M\Rightarrow  X_0AY_0\subset B 
\end{align*}
The proof is complete. $\qquad \Box$

\bigskip

One consequence of our theory is the following proposition:

\begin{proposition}\label{B10} Let $\cl M$ be a nest acting on an infinite dimensional separable Hilbert space $H.$ 
We assume that $\cl M$ is not totally atomic. Then there exists a right weak* rigged module $Y$ over $\Alg{\cl M}$ 
such that $M_l(Y)\cong \Alg{\cl N}$, where $M_l(Y)$ is the algebra of left multipliers of $Y$ and $\cl N$ 
is the nest of Volterra.
\end{proposition}
\textbf{Proof} Let $p$ be the supremum of the atoms of $\cl M.$ Our assumptions imply that  $p<I_H.$ So the nest 
$$ \cl M|_{p(H)^\bot } =\{M|_{p(H)^\bot }: M\in \cl M\}$$ is a continuous nest acting on $p(H)^\bot .$ Also, the 
map $$\sigma : \cl M\rightarrow \cl M|_{p(H)^\bot } : M\rightarrow M|_{p(H)^\bot }$$ is a continuous nest homomorphism.
All continuous nests are isomorphic with the nest of Volterra, \cite{dav}. Therefore there  exists a nest isomorphism 
$\theta : \cl M|_{p(H)^\bot } \rightarrow \cl N.$ The map $\tau =\theta\circ  \sigma: \cl M\rightarrow \cl N $ is a continuous 
surjective nest homomorphism. Put $Y=Op(\tau ).$ Theorem \ref{B1} implies that  $M_l(Y)\cong \Alg{\cl N}.$ $ \qquad \Box$

\section{Right weak* rigged modules over nest algebras}

In this section we prove that a right weak* rigged module over a nest algebra is also a left nest algebra 
module, hence  a nest algebra bimodule. We also prove that a right weak* rigged module over a nest algebra 
is $Op(\phi )$  for some continuous nest homomorphism $\phi .$ We 
need the following lemma:

\begin{lemma}\label{D1} Let $\cl M$ be a nest acting on the Hilbert space $H$,and let $B=\Alg{\cl M}$ be the 
corresponding nest algebra and $Y\subset B(H, L)$ be a weak* closed space such that $YB\subset Y.$ We also  
assume that there exists a space $X\subset B(L, H)$ such that $XY\subset B$ and $I_L\in \wsp{YX}.$ 
Then $Y$ is a nest algebra bimodule, i.e., there 
 exists a nest algebra $A\subset B(L)$ such that $AY\subset Y.$  
\end{lemma}
\textbf{Proof} Write $\theta : \cl M\rightarrow pr(B(L))$ for the map which sends $M\in \cl M$ 
to the projection onto $\overline{YM(H)}$, and denote by $A$ 
the algebra $\Alg{\theta (\cl M)}$, and put $\hat{Y}=Op(\theta ).$ Clearly $Y\subset \hat{Y}.$ 
We denote by $\hat{X}$ the space
$$\{x\in B(L, H):\;\; M^\bot x\theta (M)=0\;\;\forall \;\;M\;\in \;\cl M \}. $$
If $x\in X$ and $y\in Y$, then since $xy\in B$ we have for all $M\in \cl M$ that
$$M^\bot xyM=0\Rightarrow M^\bot x\theta (M)=0\Rightarrow x\in \hat{X}. $$
We proved $X\subset \hat{X}.$ Also for all $M\in \cl M$ 
$$ M^\bot \hat{X}\hat{Y}M =M^\bot \hat{X}\theta (M)\hat{Y}M =0\Rightarrow \hat{X}\hat{Y}\subset B. $$
Since $A\hat{Y}\subset \hat{Y} $, it suffices to prove 
$\hat{Y}= Y.$ Indeed, 
\begin{align*} \hat{Y}\;\subset \;& \wsp{YX}\hat{Y}\;\subset\; \wsp{YX\hat{Y}} \\
\subset \; &\; \wsp{Y\hat{X}\hat{Y}}\; \subset \;\wsp{YB}\;\subset\; Y. \qquad \Box
\end{align*} 

\medskip

\begin{theorem}\label{D2} Let $\cl M$ be a nest acting on the Hilbert space $H$, let $B=\Alg{\cl M}$ be the 
corresponding nest algebra, and let $Y$ be a right weak* rigged module over  $B.$ 
Then there exists a Hilbert space $K,$ a completely isometric weak* continuous $B-$module map $\Phi : Y\rightarrow  
B(H, K)$, and a nest algebra $A\subset B(K)$, such that $A\Phi (Y)B\subset \Phi (Y).$ 
\end{theorem}
\textbf{Proof} We recall the following facts from \cite{bk}. The space $K=Y\otimes ^{\sigma h}_B H$ with its 
norm is a Hilbert space. The map 
$\Phi : Y\rightarrow B(H, K)$ given by $$\Phi (y)(h)=y\otimes _Bh$$ is a completely isometric, weak* continuous $B-$module 
map. We recall the space $ \stackrel{\sim }{X } = {w^*}CB(Y, B)_B$ from Section  2. The map 
$$\Psi : \stackrel{\sim }{X } \rightarrow B(K,H), \;\;\; \Psi (u)(y\otimes _Bh)=u(y)(h)$$ is also a complete isometry 
and weak* continuous. Finally $$I_K\in \wsp{\Phi(Y) \Psi(\stackrel{\sim }{X } ) }. $$ 
 We can easily check that $$\Psi(\stackrel{\sim }{X } ) \Phi (Y)\subset B$$ and $$\Phi(Y)B\subset  \Phi(Y). $$
 So $\Phi (Y)$ satisfies the assumptions of  Lemma \ref{D1}. Therefore $\Phi (Y)$ is an $A-B$ bimodule 
for a nest algebra $A$ acting on the Hilbert space $K. \qquad \Box$

\begin{theorem}\label{B4} Let $\cl M$ be a nest acting on the Hilbert space $H$,     
let $B=\Alg{\cl M}$ be the corresponding nest algebra, and let $Y$ be a right dual operator $B-$module.  
Then the following are equivalent:

(i) The space $Y$ is a weak* rigged module over $B.$

(ii) There exists a Hilbert space $K,$ a completely isometric weak* continuous $B-$module map $\Phi : Y\rightarrow  
B(H, K)$, and a continuous nest homomorpism $\theta : \cl M\rightarrow pr(B(K))$, such that $\Phi (Y)=Op(\theta ).$     
\end{theorem}
 \textbf{Proof} 

(ii)$\Rightarrow$ (i):

By Theorem \ref{B1}, $\Phi (Y)$ is a right weak* rigged module over $B.$ Since $\Phi: Y\rightarrow  \Phi(Y) $ 
is a completely isometric weak* continuous $B-$module map, $Y$ is a right weak* rigged module over $B.$ 

(i)$\Rightarrow$ (ii):

Using  Theorem \ref{D2}, we may assume that $Y\subset B(H, L)$ for a Hilbert space $L$ and there exists a nest $\cl N$ acting on $L$ such that 
$AYB\subset Y,$ where $A=\Alg{\cl N}.$ We recall the Hilbert space $K,$ the space $\stackrel{\sim }{X } $ 
and the maps $\Phi, \Psi $ from the proof of Theorem \ref{D2}.

There exist nets 
$$ (y_t)_t\subset R^{fin}_\infty (Y) ,\;\;\; (u_t)_t\subset C^{fin}_\infty (\stackrel{\sim }{X } ) $$ 
such that $w^*-\lim_t \Phi(y_t) \Psi (u_t)=I_K.$

We define the map $$\theta : \cl M\rightarrow pr(B(K)): \theta (M)=\overline{\Phi (Y)M(H)}$$ 
which is left continuous, \cite{ep}. As in Theorem \ref{B2}, we can prove that it is right continuous.

Let $\phi $ be the restriction of $\Map{Y}$ to $\cl M.$ By Theorem \ref{A1}, $\phi (\cl M)\subset \cl N$ and $Y=Op(\phi ).$

Put $\Omega =Op(\theta ).$ It suffices to prove that $\Omega =\Phi (Y).$ Let $T$ be the operator from $K$ to 
$L$ which sends $y\otimes _Bh $ to $y(h).$ Observe that $T\Phi (y)=y$ for all $y\in Y$ and $T\theta (M)(K)\subset \phi 
(M)(H)$ for all $M\in \cl M.$ If $S\in \Omega $ and $M\in \cl M$, then
$$TSM=T\theta (M)SM=\phi (M)T\theta (M)SM =\phi (M)TSM,$$ so $TS\in Y.$ In the sequel, we can see that
$$\Phi (TS)=\lim_t\Phi(T \Phi(y_t) \Psi(u_t)S)=\lim_t\Phi (y_t\Psi (u_t)S). $$
Since $\Psi (u)S\in B$ for all $u\in \stackrel{\sim }{X } $ and $\Phi $ is a $B-$module map,
$$\Phi (TS)=\lim_t\Phi(y_t) \Psi (u_t)S=S. $$ 
We proved $\Omega \subset \Phi (Y).$ On the other hand 
$$\Phi (y)M(h)=y\otimes _BM(h)\Rightarrow \theta (M)^\bot \Phi (y)M=0 , \;\;\forall M\in \cl M. $$ So 
$\Omega \supset \Phi (Y).$ $\qquad \Box$ 

\begin{remark}\label{B5} \em{ Let $H, L$ be Hilbert spaces, $B\subset B(H)$ be a nest algebra, and 
$Y\subset B(H, L) $ be a right weak* rigged module over $B.$ We also assume that 
 $K, \Phi, \theta $ are as in the above theorem. Then there exists, see the proof, a 
contraction $T: K\rightarrow L$ such that $T\Phi (y)=y$ for all $y\in Y.$ Also we may consider that 
$\theta=\mathrm{Map} (\Phi(Y)). $} 
\end{remark}

The following two propositions can be used to construct modules over nest algebras which are not weak* rigged modules.

\begin{proposition}\label{B11} Let $\cl M$ be a continuous nest  acting on the Hilbert space $H$, $B$ be the 
 corresponding nest algebra, and $Y$ be a right weak* rigged module over $B.$ Then the algebra of left multipliers of $Y$, 
$M_l(Y),$ is isomorphic with an algebra $\Alg{\cl N}$ where $\cl N$ is a continuous nest. 
\end{proposition} 
\textbf{Proof} From Theorem \ref{B4}, there exists a weak* continuous completely isometric $B-$module map 
$\Phi : Y\rightarrow B(H, K),$ a nest algebra $A=\Alg{\cl N}\subset B(K)$, and a continuous surjective nest 
homomorphism $\theta : \cl M\rightarrow \cl N$, such that $\Phi(Y)=Op( \theta). $ By Theorem \ref{B1}, $ A\cong M_l(Y).$ 
Clearly $\cl N$ is a continuous nest.  $\qquad \Box$ 

\begin{proposition}\label{B12} 
Let $\cl M$ (resp. $\cl N$) be a nest   acting on the Hilbert space $H$ (resp. $K$), let $B,$ (resp. $A$) be the 
 corresponding nest algebra, and let $Y$ be a weak* closed $A-B$ bimodule, such that $\overline{Y(H)}=K.$ We 
assume that $\cl M$ is a continuous nest and $\cl N$ has at least one atom. Then $Y$ is not a right weak* rigged module over $B.$ 
\end{proposition} 
\textbf{Proof} Let  $q$ be an atom of $\cl N.$ It follows that the space $qB(K)q$ is a subset of the algebra $A.$ 
The space $qY$ is a subspace  of $B(H, q(K))$ and satisfies 
$$B(q(K))qYB\subset qY. $$
So $qY$ is a nest algebra bimodule over the algebras $B(q(K)),$ corresponding to the nest $ \{0_K, q\} ,$ 
and $B.$ If $$\theta =\Map{qY}: \cl M\rightarrow \{0_K, q\} $$ by Theorem \ref{A1} 
$$qY=\{y\in B(H, q(K)): \theta (M)^\bot yM=0,\;\;\forall \;M\;\;\in \;\;\cl M\}. $$
 We define the projection $$M_0=\vee \{M\in \cl M:\;\; \theta (M)=0\}. $$ Since $\theta $ is 
left continouous, $\theta (M_0)=0.$ If $M>M_0,$ then $\theta (M)>0.$ Thus $\theta (M)=q.$ 
 Therefore $$qY=\{y\in B(H, q(K)): yM_0=0\}=B(H, q(K))M_0^\bot . $$
Assume now that $Y$ is a weak* rigged module over $B.$ Using for example Definition \ref{themelion}, we 
can verify that $qY$ is a weak* rigged module over $B.$ In this case, 
$$M_l(qY)=B(q(K))=\Alg{\{0_K, q\}}. $$ This contradicts Proposition \ref{B11}. $\qquad \Box$

\section{Examples}

\begin{example}\label{C1} If $ B=\Alg{\cl M} \subset B(H)$ is a nest algebra corresponding to a finite nest $\cl M,$   
$A=\Alg{\cl N}$ is a nest algebra, and $Y$ is a  weak* closed $B-A$ bimodule, then $Y$ is a right weak* rigged module over $B.$ 
Indeed, the map $\phi : \cl M\rightarrow \cl N$ sending every $M\in \cl M$ to the projection onto $\overline{YM(H)}$
 is  continuous. Use now Theorems \ref{A1}, \ref{B1}.
\end{example}

\begin{example}\label{C2} Let $B=\Alg{\cl M}, A=\Alg{\cl N}$ be  nest algebras and $\phi : \cl M\rightarrow \cl N$ 
be  a nest isomorphism. Since $\phi $ is  continuous, the space $Y=Op(\phi )$ is a right weak* rigged module over $B.$ 
\end{example}

\begin{example}\label{C3} Let $H$ be an infinite dimensional separable Hilbert space and $P_n$ be a 
strictly increasing sequence 
of projections such that $\vee _{n}P_n=I_H.$ 
The set $\cl M=\{0_H, P_n, n\in \mathbb{N}, I_H\}$ is a nest. Suppose that $A=\Alg{\cl N}$ is another nest algebra 
and $Y$ is a weak* closed $A-\Alg{\cl M} $ bimodule. Then $Y$ is a right weak* rigged module over $\Alg{\cl M} .$ Indeed, the map 
$$\phi =\Map{Y}: \cl M\rightarrow \cl N$$ is continuous in every $P_n$ and left continuous in $I_H,$ by Theorem \ref{A1}. 
Therefore $\phi $ is continuous. 
\end{example}

\begin{example}\label{C5} Let $B$ be a nest algebra and $s$ be an invertible operator.  
By Theorem \ref{B6}, $Y=sB$ is a right weak* rigged module over $B.$
\end{example}

\begin{example}\label{C7} 
Let $\{e_n: n\in \mathbb{N}\}$ be an orthonormal basis  of the Hilbert space $H$ and 
 $p_n$ be the projection onto the space generated by the vectors $\{e_1, e_2, ...,e_n\}$ for all $n.$ Let 
$\cl M$ be the nest $\{0_H, p_n, I_H, \;\;n\;\in \;\mathbb{N}\}, $ $B$ be the algebra $\Alg{\cl M}$, and 
$x, y$ be the shift operators given by $$y(e_1)=0, y(e_n)=e_{n-1}, \;\;\forall n\geq 2,$$$$ x(e_n)=e_{n+1}, \;\;\forall n\geq 
1. $$ Since $x^m y^m$ is the projection onto $p_m(H)^\bot $ which belongs to $B$ and 
$y^mx^m=I_H$ for all $m, $  Theorem \ref{B6} then implies that $Y=\overline{y^mB}^{w^*}$ is a right weak* rigged module over $B$ 
 for all $m.$  
\end{example}

\begin{example}\label{C20}Let $B$ be a nest algebra and $\cl T$ be a ternary ring of operators such that $\cl T^* \cl T
\subset B.$ Since   $\wsp{\cl T \cl T^*}$ is a unital algebra by Theorem \ref{B6} the space 
$\wsp{\cl T B}$ is a weak* rigged module over $B.$ Modules of this type are called {\em projective weak* rigged modules}
 and they have 
special properties \cite{bkraus}.  
\end{example}

The next example is an example of a nest algebra bimodule which is not a weak* rigged module.

\begin{example}\label{C8}
Let $\cl M$ be a continuous nest acting on $H, B$ be the algebra $\Alg{\cl M}$, and $Y=B(H, K)$ where $K$ is another 
Hilbert space. The space $Y$ is a nest algebra bimodule over $B(K)$ and $B.$ But it is not a weak* rigged module over 
$B$ because $\cl M$ is a continuous nest and the nest $\{0_K, I_K\}$ corresponding to $B(K)$ is totally atomic. 
 (Use Proposition \ref{B12}).   
\end{example}

This example is also an example of an $A-B$ bimodule $Y$ where both $A$ and $B$ are continuous nest algebras 
and $Y$ is not a weak* rigged module over $B:$ Choose infinite dimensional Hilbert spaces $K$ and $H$, 
and continuous nest algebras $B\subset B(K)$ and $A\subset B(L).$ In the same way we can prove that $Y=B(H, K)$  
is not a weak* rigged module over $B.$

\begin{example}\label{C9}In Example \ref{C3}, $\cl M$ is a totally atomic nest and $Y$ is a right weak* rigged 
module over $\Alg{\cl M}.$ As we can easily prove, $M_l(Y)$ 
 is isomorphic to a nest algebra $\Alg{\cl N}$ 
where $\cl N$ is  not a continuous nest. Here, 
we prove that there is a totally atomic nest $\cl M$ and a right weak* rigged module $Y$ over $B=\Alg{\cl M}$ such that 
$M_l(Y)\cong \Alg{\cl N_v}$ where  
$\cl N_v$ is the nest of Volterra:

 There is a totally atomic nest $\cl M$ (the nest of Cantor) and an isomorphic nest 
$\cl M_2$ which is not totally atomic,  \cite[Example 7.19]{dav}. Suppose that $\rho : \cl M\rightarrow \cl M_2$ is the above 
nest isomorphism. Let $\tau : \cl M_2\rightarrow \cl N_v$ be the surjective continuous nest homomorphism 
described in Proposition \ref{B10}. If $Y=Op(\tau \circ \rho ),$ Theorem \ref{B1} implies that $M_l(Y)\cong \Alg{\cl N_v}.$ 
\end{example}

\section{Spatially embedding algebras}

In this section we investigate  a weaker definition than  Definition \ref{A2}. The objects of the theory  are the 
reflexive algebras and more specifically the CSL algebras. A nest algebra is a special type of a CSL algebra. 
 See \cite{dav} or in \cite{st} for the definition 
of reflexivity and the notions of CSL, the  CSL algebra and a  reflexive algebra.  

\begin{definition}\label{E1} Let $A$ (resp. $B$ ) be a weak* closed algebra acting on $H$ (resp. $K$). We say 
that $A$ is \textbf{spatially embedded} in $B$ if there exist a $B-A$ bimodule $X\subset B(H,K)$ and an $A-B$ bimodule 
$Y\subset B(K,H)$ such that
 
(i) $XY\subset B$,

(ii) $A=\wsp{YX}.$

Moreover, if  $B=\wsp{XY}$, we call $A$ and $B$ \textbf{spatially Morita equivalent}.
\end{definition}

\begin{remark}\label{E2}Let $\cl M$ and $ \cl N$ be nests corresponding to the algebras $B$ and $ A$ and 
let $\phi : \cl M\rightarrow \cl N$ be  a continuous surjective nest homomorphism. It follows that  $A$ is embedded spatially in $B$ 
(see Theorem \ref{B2}). On the other hand, if $B=\Alg{\cl M}$ is a nest algebra and $Y$ is a right 
weak* rigged module over $B$, there is a normal completely isometric representation $\Phi $ of $Y$ 
and a nest algebra $A=\Alg{\cl N}$ such that $\Phi (Y)$ is an $A-B$ bimodule and $A$ is spatially embedded in $B$ 
(see Theorems \ref{D2} and \ref{B4}). In this case there exists a continuous nest homomorphism from $\cl M$ 
onto $\cl N.$
\end{remark}

\begin{proposition}\label{E3} Let $A, B, X, Y$ be as in Definition \ref{E1}. Put  
$$\Omega_2=M_2(B),\;\;\;\;  \Omega_1= \left(\begin{array}{clr} B & X \\ Y & A\end{array}\right)  . $$
The algebras   $\Omega_1$ and $ \Omega_2 $ are spatially Morita equivalent.

\end{proposition}
\textbf{Proof} We define spaces
$$\hat{Y}=\left(\begin{array}{clr} B & B \\ Y & Y\end{array}\right) ,\;\;\; \hat{X}=
\left(\begin{array}{clr} B & X \\ B & X\end{array}\right).  $$
We can easily check that $\hat{Y}$ is an $\Omega_1- \Omega_2 $ bimodule,  $\hat{X}$ is an $\Omega_2- \Omega_1 $ 
bimodule, and $$ \Omega _2=\wsp{\hat{X}\hat{Y}}, \;\;\;\; \Omega _1=\wsp{\hat{Y}\hat{X}}. \qquad \Box$$

\begin{corollary}\label{E4} Let $A, B, X, Y$ be as in Definition \ref{E1}. If $B$ is a reflexive 
algebra, then $X, Y$, and $A$ are  reflexive spaces.
\end{corollary}    
\textbf{Proof} If $B$ is reflexive,  $M_2(B)$ is reflexive. It follows from \cite[remark 4.2]{ele} 
that the algebra $\Omega _1$ defined in Proposition \ref{E3} is reflexive. Therefore 
$X, Y$, and $A$ are reflexive. $\qquad \Box$

\medskip

In the rest of this section, if $Z$ is a set of operators, then $\cl K(Z)$ denotes its subset 
of compact operators, $\cl F(Z)$ denotes its subset of finite rank operators, and  $\cl R(Z)$ denotes its subset of rank one operators.

\begin{theorem}\label{compact} Let $A, B, X, Y$ be as in Definition \ref{E1}. 

(i) If $B=\overline{\cl K(B)}^{w^*}$, then  $A=\overline{\cl K(A)}^{w^*}.$

(ii) If $B=\overline{\cl F(B)}^{w^*}$, then  $A=\overline{\cl F(A)}^{w^*}.$

(iii) If $B=\wsp{\cl R(B)}$, then  $A=\wsp{\cl R(A)}.$

\end{theorem}
\textbf{Proof} Suppose that $B=\overline{\cl K(B)}^{w^*}.$ Since $Y$ is a left $B$-module, 
$$Y=\wsp{YB}=\wsp{Y\cl K(B)}=\overline{\cl K(Y)}^{w^*}. $$ Therefore 
$$A=\wsp{YX}=\wsp{\cl K(Y)X}=\overline{\cl K(A)}^{w^*}. $$
The proofs of (ii) and (iii) are  similar. $\qquad \Box$

\begin{theorem}\label{compop} Let $A, B, X, Y$ be as in Definition \ref{E1}. Assume that $A$ is a unital algebra. 

(i) If $s\in \cl K(A)$ is a nonzero operator, then $B$ contains a nonzero compact operator. 

(ii) If $s\in \cl F(A)$  is a nonzero operator, then $B$ contains a nonzero finite rank operator. 

(iii) If $s\in \cl R(A)$,  then $B$ contains a rank one operator. 
\end{theorem}
\textbf{Proof} Let $s\in \cl K(A)$ be a nonzero operator. There exist $x_0\in X, y_0\in Y$ such that $x_0sy_0\neq
0.$ Indeed, if $XsY=0$ then $YXsYX=0.$ By assumption, $\wsp{YX}$ contains the identity operator. So $s=0.$ 
This is a contradiction. 

Observe that $0 \neq x_0sy_0\in \cl K(B),$ so statement (i) holds. The proofs of statements (ii), (iii) are similar. 
$\qquad \Box$

\begin{theorem}\label{E5} Let $\cl L_1, \cl L_2$ be CSLs acting on the Hilbert spaces 
$K, H$ respectively and $A, B$ be the corresponding CSL algebras. The following are equivalent:

(i) $A$ is embedded spatially in $B.$

(ii) There exists a continuous surjective lattice homomorphism $\phi : \cl L_2\rightarrow \cl L_1.$
\end{theorem}
\textbf{Proof}

(i)$\Rightarrow $(ii):

Suppose that there exist spaces $X$ and $Y$ satisfying (i) and (ii) of the definition \ref{E1}. We define the spaces 
$\hat{X}, \hat{Y}$ and the algebras $\Omega_1,  \Omega_2 $ as in Proposition \ref{E3}. In the sequel if $C$ is a unital 
algebra, we denote by $\Lat{C}$ the lattice of invariant projections of  $C.$ If 
$\sigma =\Map{\hat{Y}},$ Theorem 4.1 in \cite{ele} implies that $\sigma $ is a lattice isomorphism from $\Lat{\Omega _2}$ 
 onto $\Lat{\Omega _1}.$ 

If $P\in \cl L_2$ we denote  by $\theta (P)$ the projection onto $\overline{YP(H)}$, which clearly 
belongs to $\cl L_1.$ We can easily check that 
$$\sigma((P\oplus P))= P\oplus \theta (P)$$ for all $P\in \cl L_2.$ So the map 
$\theta : \cl L_2\rightarrow \cl L_1$ is a continuous CSL homomorphism. It remains to show that $\theta $ is surjective. If 
$Q\in \cl L_1$, then $\overline{AQ(K)}=Q(K).$ The projection onto the closure of the linear span of the set
$$\{\Omega _1 \left(\begin{array}{clr}0& 0\\ 0 & Q\end{array}\right)
\left(\begin{array}{clr}\xi \\ \eta \end{array}\right) : \;\;\xi \;\in\; H, \;\;\;\eta \;\in K \}$$
belongs to 
$$\Lat{\Omega _1}=\{(P\oplus \theta (P)): \;\;P\;\;\in\;\;\cl L_2 \}. $$
 So it is equal to a projection  $R=(P_0\oplus \theta (P_0))$ for $P_0\in \cl L_2.$

 Observe  that $R$ is the projection onto the space generated by the vectors of the form
$$ \left(\begin{array}{clr}0& xQ\\ 0 & aQ\end{array}\right)\left(\begin{array}{clr}\xi \\ \eta \end{array}\right) 
: \;\;\xi \;\in\; H, \;\;\; \eta \;\in K , x \in X, \;\; a\in A . $$
So $\theta (P_0)$ is the projection onto $\overline{AQ(K)}$ which is $Q.$ Therefore $\theta $ is surjective. 

\medskip

(ii)$\Rightarrow$ (i):

Let $\phi : \cl L_2\rightarrow \cl L_1$ be a continuous surjective lattice homomorphism. We define the CSLs 
$$ \cl N_2=\{P\oplus P: \;\;P\;\in \;\;\cl L_2\},  \;\;\;\cl N_1=\{P\oplus \phi (P): \;\;P\;\in \;\;\cl L_2\}.  $$
 The map  $\rho : \cl N_2\rightarrow \cl N_1$ sending 
every $P\oplus P$ to $P\oplus \phi (P)$ is a CSL isomorphism. We define the spaces

\begin{align*}\hat{Y}=&\{y\in B(K\oplus K, K\oplus H): \rho (L)^\bot yL=0, \;\;\forall L\in \cl N_2 \},\\  
\hat{X}=&\{x\in B(K\oplus H, K\oplus K): L^\bot x\rho (L)=0, \;\;\forall L\in \cl N_2 \}, \\  
Y=&\{y\in B(K,H): \phi (P)^\bot yP=0, \;\;\forall P\in \cl L_2 \}, \\  
X=&\{x\in B(H,K): P^\bot x\phi (P)=0, \;\;\forall P\in \cl L_2 \}. \end{align*}

We can easily verify that 
 $$\hat{Y}= \left(\begin{array}{clr} B & B \\ Y & Y\end{array}\right) ,\;\;\; \hat{X}=
\left(\begin{array}{clr} B & X \\ B & X\end{array}\right),\;\;\;
\Alg{\cl N_1}=   \left(\begin{array}{clr} B & X \\ Y & A\end{array}\right).  $$
By Theorem 4.3 in \cite{ele}, we conclude that $$\Alg{\cl N_1}=\wsp{\hat{Y}\hat{X}},$$ therefore $A=\wsp{YX}.$ 
The proof is  complete. $\qquad \Box$

\medskip
There is no known lattice condition corresponding to the weak* density of the compact operators in 
 a CSL algebra. But comparing Theorems \ref{compact} and \ref{E5}, we have the following corollary:

\begin{corollary}\label{density} Let $\cl L_1, \cl L_2$ be CSLs and $\phi : \cl L_2\rightarrow \cl L_1$ 
be a continuous surjective lattice homomorphism. 

(i)If $ \mathrm{Alg}(\cl L_2) =\overline{\cl K(\mathrm{Alg}(\cl L_2) )}^{w^*}$,
then $ \mathrm{Alg}(\cl L_1) =\overline{\cl K( \mathrm{Alg}(\cl L_1) )}^{w^*}.$

(ii)If $ \mathrm{Alg}(\cl L_2) =\overline{\cl F(\mathrm{Alg}(\cl L_2) )}^{w^*}$,
then $ \mathrm{Alg}(\cl L_1) =\overline{\cl F( \mathrm{Alg}(\cl L_1) )}^{w^*}.$
\end{corollary}

\medskip
\medskip
Similarly Theorems \ref{compop} and \ref{E5} imply

\begin{corollary}Let $\cl L_1, \cl L_2$ be CSLs and $\phi : \cl L_2\rightarrow \cl L_1$ 
be a continuous surjective lattice homomorphism. If $s\in \mathrm{Alg}(\cl L_1) $ is, respetively, a nonzero compact operator, a 
nonzero finite rank operator, a rank one operator, then $\mathrm{Alg}(\cl L_2) $ contains, respectively, a  
nonzero compact operator, a nonzero finite rank operator, a rank one operator.  
\end{corollary}

\medskip

\medskip

If $U$ is a reflexive and separably acting bimodule over maximal abelian selfadjoint algebras ({\em masa})s, there exists a 
smallest weak* closed masa bimodule $U_{min}$ whose reflexive hull is $U.$ In the special case $U=U_{min}$, we call $U$ 
{\em synthetic}. A CSL $\cl L$ is called {\em synthetic} if the algebra $\Alg{\cl L}$ is synthetic \cite{dav}. We 
present the following relevant result:

\begin{theorem}\label{E6} Let $\cl L_1$ and $ \cl L_2$ be separably acting CSLs, corresponding to the algebras $B$ and $ A$, and suppose $\phi : \cl L_2\rightarrow \cl L_1$ is 
a continuous surjective lattice homomorphism. If $\cl L_2$ is synthetic, then $\cl L_1$ is synthetic.
\end{theorem}  
\textbf{Proof} We define the CSLs

$$ \cl N_2=\{P\oplus P: \;\;P\;\in \;\;\cl L_2\},  \;\;\;\cl N_1=\{P\oplus \phi (P): \;\;P\;\in \;\;\cl L_2\}.  $$
The map $$\rho : \cl N_2\rightarrow \cl N_1,\;\;\; P\oplus P\rightarrow P\oplus \phi (P)$$ is a CSL isomorphism.
We have that  
\begin{equation} \label{csl1} \Alg{\cl N_1}=   \left(\begin{array}{clr} B & Y \\ X & A\end{array}\right), \end{equation}
 where 
$X$ and $Y$ are defined as in Theorem \ref{E5}. As in the proof of Theorem 4.7 in \cite{ele}, we can prove that 
 \begin{equation} \label{csl2} \Alg{\cl N_1}_{min}=   \left(\begin{array}{clr} B_{min} & Y_{min} 
\\ X_{min} & A_{min}\end{array}\right).\end{equation}
 
If $\cl L_2$ is synthetic then $\cl N_2$ is synthetic. By Theorem 4.7 in \cite{ele}, $\cl N_1$ is 
synthetic, so by (\ref{csl1}) and  (\ref{csl2}),  $A=A_{min}. \qquad \Box$

\end{document}